\documentclass[12pt,reqno]{amsart}
\usepackage{amssymb}
\makeatletter
\def\proof{\@ifnextchar[{\@oproof}{\@nproof}}
\def\@oproof[#1][#2]{\trivlist\item[\hskip\labelsep
\textit{#2 Proof of\ #1.}~]\ignorespaces}
\def\@nproof{\trivlist\item[\hskip\labelsep\textit{Proof.}~]\ignorespaces}

\makeatother

\setbox0=\hbox{$+$}
\newdimen\plusheight
\plusheight=\ht0
\def\+{\;\lower\plusheight\hbox{$+$}\;}

\setbox0=\hbox{$-$}
\newdimen\minusheight
\minusheight=\ht0
\def\-{\;\lower\minusheight\hbox{$-$}\;}

\setbox0=\hbox{$\cdots$}
\newdimen\cdotsheight
\cdotsheight=\plusheight%\ht0
\def\cds{\lower\cdotsheight\hbox{$\cdots$}}

\setbox0=\hbox{$+$}

\renewcommand{\(}{\left\(}
\renewcommand{\)}{\right\)}
\renewcommand{\[}{\left[}

\numberwithin{equation}{section}
 \theoremstyle{plain}
\newtheorem{thm}{Theorem}[section]

\newtheorem{rem}[thm]{Remark}

\begin{document}
\title[A Simple Number Theoretic Result] {A Simple Number Theoretic Result}
\author{Manjil P.~Saikia}
\address{Department of Mathematical Sciences, Tezpur University, Napaam, Sonitpur, Pin-784028, India}
\email{manjil\_msi09@agnee.tezu.ernet.in, manjil.saikia@gmail.com}
\author{Jure Vogrinc}
\address{Faculty of Mathematics and Physics, University of Ljubljana, Jadranska ul.~19, 1000, Ljubljana, Slovenia}
\email{jure.vogrinc@gmail.com}

%\maketitle

\vspace*{0.5in}
\begin{center}
{\bf A Simple Number Theoretic Result}\\[5mm]
{\footnotesize  MANJIL P.~SAIKIA\footnote{Corresponding Author: manjil.saikia@gmail.com} and JURE VOGRINC}\\[3mm]
\end{center}

\vskip 5mm \noindent{\footnotesize{\bf Abstract.} We derive an interesting congruence relation motivated by an Indian Olympiad problem. We give three different proofs of the theorem and mention a few interesting related results.}

\vskip 3mm


\vskip 3mm

\noindent{\footnotesize 2000 Mathematical Reviews Classification
Numbers: 11A07, 11A41, 11A51.}

\section{An Interesting Olympiad Problem}

\begin{thm}
$7$ divides $\binom{n}{7}-\lfloor\frac{n}{7}\rfloor$, $\forall n \in \mathbb{N}$.
\end{thm}

The above appeared as a problem in the Regional Mathematical Olympiad, India in 2005. Later, in 2007, a similar type of problem was set in the undergraduate admission test of Chennai Mathematical Institute, a premier research institute of India where $7$ was replaced by $3$.

In late 2008, the first author posted the following theorem as a question in an internet forum called Mathlinks, \cite{ml}.

\begin{thm}
If $p$ is any prime then, $p$ divides $\binom{n}{p}-\lfloor\frac{n}{p}\rfloor$, $\forall n \in \mathbb{N}$.
\end{thm}

The second author replied to the post and proved the above result using Wilson's theorem. However, his proof was not entirely correct and together the authors managed to correct the argument.

Then we wondered what would be the case if $p$ is not a prime.

In early 2009, we showed that if $p$ is a composite number of the form $q^x.k$, where $q$ is a prime and $gcd(q,k)=1$, then the above statement is not true. Note here that $x$ and $k$ cannot simultaneously be equal to $1$.

Thus \textbf{Theorem 1.2} became the basis of the following theorem.

\section{{Main Result}}

\begin{thm}
A natural number $p>1$ is a prime if and only if $\binom{n}{p}-\lfloor\frac{n}{p}\rfloor$ is divisible by $p$ for every non-negative $n$, where $\binom{n}{p}$ is the number of different ways in which we can choose $p$ out of $n$ elements and $\lfloor x \rfloor$ is the greatest integer not exceeding the real number $x$.
\end{thm}

We have found three different proofs of the above result, two purely number theoretic and one via a combinatorial argument.

We also state and prove a famous theorem in Number Theory called Lucas' Thoerem,

\begin{thm}
\textbf{(\cite{el},E.~Lucas, 1878)}Let $p$ be a prime and $m$ and $n$ be two integers considered in the following way,\\

$$m=a_kp^k+a_{k-1}p^{k-1}+\ldots+a_1p+a_0$$
$$n=b_lp^l+b_{l-1}p^{l-1}+\ldots+b_1p+b_0,$$

where all $a_i$ and $b_j$ are non-negative integers less than $p$. Then,

$$\binom{m}{n}=\binom{a_kp^k+a_{k-1}p^{k-1}+\ldots+a_1p+a}{b_lp^l+b_{l-1}p^{l-1}+\ldots+b_1p+b_0} \equiv \Pi_{i=0}^{max(k,l)}\binom{a_i}{b_i}.$$
\end{thm}

There are numerous proofs of the above theorem, the first being given by Lucas himself. We give a slightly modified version of a proof given by N.~J.~Fine, \cite{njf}.

\begin{proof}
All we need to prove is that if $p$ is a prime than $\binom{ap+b}{cp+d}=\binom{a}{c}\binom{b}{d}$, for non negative $a,b,c,d$ where $b,d<p$. Once we have this, then we can use induction on it to get the desired result.

We have, for every integer $k$, such that $0< k <p$ we have $p\mid \binom{p}{k}$, so we can conclude that for every integer $x$ we have,

$$(1+x)^p\equiv(1+x^p)~(mod~p).$$

So, from the above result we can easily see that,

$$(1+x)^{ap+b}=((1+x)^p)^a)(1+x)^b\equiv(1+x^p)^a(1+x)^b~(mod~p).$$

Comparing the polynominal coefficients of $x^{cp+d}$ in both congruent polynomials, we get,

$$\binom{ap+b}{cp+d}=\binom{a}{c}\binom{b}{d}.$$

This completes the proof.
\end{proof}

Using Lucas' Theorem, we can easily prove our result, if we consider $n=ap+b$, where $p$ is a prime and $b<p$.

We now give a combinatorial proof of \textbf{Theorem 2.1},

\begin{proof}
We suppose $n=pq+r$, and let us have $n$ compartments with partitions after every $p$ compartments. In all there are $q+1$ partitions where the last partition has $r$ compartments. It is assumed that each compartment is labelled and distinct. Obviously, $\binom{n}{p}$ is the number of ways to select $p$ compartments out of $n$ compartments.

The total number of selections thus can be divided into two sets:
\begin{itemize}
\item Compartments of only one partition are completely selected.
\item Compartments of only one partition are not selected.
\end{itemize}

It is easy to see that there are $q=\lfloor\frac{n}{p}\rfloor$ selections of the first set.

In the second set we subdivide the selections on the basis of number of selections in each partition. Let $a_i$ compartments be selected in the $i-th$ partition. For each subdivision, the number of selections are $\binom{r}{a_{q+1}}\Pi_{i=1}^q\binom{p}{a_i}$.

Since $a_i<p$ and at least one $a_i(1\leq i\leq q)$ will be non zero and $\binom{p}{k}$ is divisible by $p$ for $0<k<p$.

Hence, each of the subdivisions of the second set will be divisible by $p$.

This completes the proof.

\end{proof}

We give now the final proof of $\textbf{Theorem 2.1}$ using purely number theoretic tools.

\begin{proof}
First assume that $p$ is prime. Now we consider $n$ as $n=ap+b$ where $a$ is a non-negative integer and $b$ an integer $0\leq b<p$. Obviously,

$$\lfloor \frac{n}{p}\rfloor=\lfloor \frac{ap+b}{p}\rfloor\equiv a~(mod~p).$$

Now let us calculate $\binom{n}{p}~(mod~p)$.

$$\binom{n}{p}=\binom{ap+b}{p}$$
$$=\frac{(ap+b)\cdot(ap+b-1)\cdots(ap+1)\cdot ap\cdot(ap-1)\cdots(ap+b-p+1)}{p\cdot(p-1)\cdots 2\cdot1}$$
$$=\frac{a\cdot(ap+b)\cdot(ap+b-1)\cdots(ap+1)\cdot(ap-1)\cdots(ap+b-p+1)}{(p-1)\cdot(p-2)\cdots 2\cdot 1}$$

We denote this number by $X$.

We have $X\equiv c~(mod~p)$ for some $0\leq c<p$. Consequently taking modulo $p$, we have
 $$c(p-1)!=X(p-1)!=a(ap+b)\cdots(ap+1)(ap-1)\cdots(ap+b-p+1)$$

All the numbers $ap+b,\dots ,ap+b+1-p$ (other than $ap$) are relatively prime to $p$ and obviously none differ more than $p$ so they make a reduced residue system modulo $p$, meaning we have mod $p$,
$$(p-1)!=(ap+b)\cdot(ap+b-1)\dots(ap+1)\cdot(ap-1)\cdot(ap+b-p+1)$$
both sides of the equation being relatively prime to $p$ so we can deduce $X\equiv c \equiv a~(mod~p)$. And finally  $\binom{n}{p}\equiv X\equiv a\equiv \lfloor \frac{n}{p}\rfloor~(mod~p)$.

To complete the other part of the theorem we must construct a counterexample for every composite number $p$. If $p$ is composite we can consider it as $q^x\cdot k$ where $q$ is some prime factor of $p$, $x$ its exponent and $k$ the part of $p$ that is relatively prime to $q$ ($x$ and $k$ cannot be simultaniously $1$ or p is prime). We can obtain a counterexample by taking $n=p+q=q^xk+q$ will make a counter example. We have: $$\binom{p+q}{p}=\binom{p+q}{q}=\frac{(q^xk+q)(q^xk+q-1)\dots (q^xk+1)}{q!}$$
Which after simplifying the fraction equals: $(q^{x-1}k+1)\frac{(q^xk+q-1)\dots (q^xk+1)}{(q-1)!}$. Similary as above we have $(q^xk+q-1)\dots (q^xk+1)=(q-1)!\neq 0$ modulo $q^x$ therefore,
$$\frac{(q^xk+q-1)\dots (q^xk+1)}{(q-1)!}\equiv 1~(mod~q^x)$$ and
 $$\binom{p+q}{p}\equiv q^{x-1}k+1~(mod~q^x).$$

On the other hand obviously,
$$\lfloor\frac{q^xk+q}{q^xk}\rfloor\equiv 0~(mod~q^x).$$

And since $q^{x-1}k+1$ can never be equal to $0$ modulo $q^x$ we see that $$\binom{p+q}{p}\neq \lfloor\frac{p+q}{p}\rfloor~(mod~q^x)$$
consequently also incongruent modulo $p=q^xk$.

\end{proof}

\begin{rem}
Here we would like to comment that by taking $q$ as the minimal prime factor of $p$ and using the same method as above we can simplify the proof even more. We can than compare $\lfloor\frac{p+q}{p}\rfloor$ and $\binom{p+q}{p}$ directly modulo $p=q^xk$ and not $q^x$.
\end{rem}

\begin{rem} In \cite{kayal}, N.~Kayal mentioned that \textbf{Theorem 2.1} can be considered to be a very naive and crude primality test.
\end{rem}

\begin{rem} Instead of looking modulo $p$, we can look at higher powers of $p$, or we can look at the $n$-th Fibbonacci prime and so on. However, initial investigations by the authors suggest that finding a congruence relation in those cases becomes more difficult.
\end{rem}

\section{{Other Interesting Results}}

It is worthwhile to mention without proof a few interesting results that may be deduced from \textbf{Lucas' Theorem} or \textbf{Theorem 2.1}.

\begin{thm}
\textbf{(\cite{rm},Romeo Me$\check{s}$trovi$\acute{c}$, 2009)}If $d,q>1$ are integers such that,

$$\binom{nd}{md}\equiv \binom{n}{m}~(mod~p)$$

for every pair of integers $n\geq m \geq 0$, then $d$ and $q$ are powers of the same prime $p$.
\end{thm}

\begin{rem}
It follows from the Lucas' Theorem, the congruence relation $\binom{np}{mp}\equiv \binom{n}{m}~(mod~p)$.
\end{rem}

\begin{rem}
We do not yet know which are the possible powers of $p^k$ and $p^l$ such that

$$\binom{np^k}{mp^k}\equiv \binom{n}{m}~(mod~p^l)$$

for all integers $n\geq m\geq 1$.
\end{rem}
\begin{thm}
\textbf{(\cite{db}, David F.~Bailey, 1990)} Let $n$ and $r$ be non-negative integers and $p\geq 5$ be a prime. Then,

$$\binom{np}{rp}\equiv \binom{n}{r}~(mod~p^3),$$

where we set $\binom{n}{r}=0$, if $n<r$.
\end{thm}

\begin{thm}
\textbf{(\cite{db}, David F.~Bailey, 1990)} Let $N,R,n$ and $r$ be non-negative integers and $p\geq 5$ be a prime. Suppose $n,r<p$, then

$$\binom{Np^3+n}{Rp^3+r}\equiv\binom{N}{R}\binom{n}{r}~(mod~p^3).$$
\end{thm}

\begin{thm}
\textbf{(\cite{tma}, Tom M.~Apostol)} If $p$ is a prime, $\alpha$ is a positive integer and $p^\alpha\mid \lfloor\frac{n}{p}\rfloor$ then, $p^\alpha \mid \binom{n}{p}$.
\end{thm}

\section{{Acknowledgements}}
The authors would like to thanks Prof.~Nayandeep Deka Baruah for reading through an earlier version of the work and pointing to us many mistakes and also to the anonymous referee for various useful suggestions/corrections.


\begin{thebibliography}{99}

\bibitem{tma}T.~M.~Apostol, \emph{An Introduction to the Analytic Theory of Numebrs}, Springer--Verlag, 1975.
\bibitem{njf} N.~J.~Fine, \emph{Binomial Coefficients modulo a prime}, Amer.~Math.~Monthly, \textbf{54} (1947), 589--592.
\bibitem{el} E.~Lucas, \emph{Th$\acute{e}$orie des Fonctions Num$\acute{e}$riques Simplement P$\acute{e}$riodiques}, Amer.~J.~Math., \textbf{1} (2), 184--196; \textbf{1} (3), 197--240; \textbf{1} (4), 289--321 (1878).
\bibitem{kayal}N.~Kayal, \emph{Primality Testing}, ICM 2010 Satellite Int.~Conf.~on Rings and Near Rings, North--Eastern Hill University, Shillong, India, 2010.
\bibitem{rm} R.~Me$\check{s}$trovi$\acute{c}$, \emph{A Note on the Congruence $\binom{nd}{md}\equiv \binom{n}{m}~(mod~q)$}, Amer.~Math.~Monthly, \textbf{116} (2009), 75--77.
\bibitem{nehu} M.~P.~Saikia, \emph{A Few Results in Number Theory}, ICM 2010 Satellite Int.~Conf.~on Rings and Near Rings, North--Eastern Hill University, Shillong, India, 2010.
\bibitem{ml}M.~P.~Saikia, J.~Vogrinc, \emph{Let's Generalize}, MathLinks Forum Discussion, 2008.
\bibitem{db} J.~Zhao, \emph{Bernoulli Numbers, Wolstenholme's Theorem, and $p^5$ Variations of Lucas' Theorem}, preprint.

\end{thebibliography}
\end{document}